\documentclass[onecolumn,notitlepage,aps,pra,10pt]{revtex4-1}
\usepackage{graphicx} 
\usepackage{amsmath}

\begin{document}
\title{Geometric considerations for energy minimization of topological links and chainmail networks}
\author{Alexander R. Klotz }
\affiliation{Department of Physics and Astronomy, California State University, Long Beach}

\begin{abstract}

Knot and link energies can be computed from sets of closed curves in three dimensional space, and each type of knot or link has a minimum energy associated with it. Here, we consider embeddings of links that locally or globally minimize the Möbius and Minimum Distance energies. By describing these energies as functions of a small number of free parameters, we can find configurations that minimize the energies with respect to these parameters. It has previous been demonstrated that such minimizers exist, but the specific embeddings have not necessarily been found. We find the geometries leading to minimal configurations of Hopf links and Borromean rings, as well as more complex structures such as chain links and chainmails. We find that scale-invariant properties of these energies can lead to ``non-physical'' minimizers, e.g. that a linear chain of Hopf links will subtend a finite length as its crossing number diverges. This incidentally allows us to derive a conjectural improved universal lower bound for the ropelength of knots and links. We also show that Japanese-style square chainmail networks are more efficient, in terms of excess energy, than square lattice ``4-in-1'' chainmail networks. 
\end{abstract}

\maketitle

\section{Introduction}

How does one produce a ``nice'' embedding of a knot in three-dimensional space? One proposed method is define an energy functional over all possible embeddings of the knot, and then to find the embedding that minimizes that energy. Knot energies typically have the following desired properties:
\begin{itemize}
    \item Each knot type has a real, finite minimum energy.
    \item The energy of a conformation is invariant under rotations, translations, and rescalings of the knot.
    \item There are a finite number of knots with minima below a certain energy, such that energy minimization can be used to classify knots.
    \item The energy diverges as parts of the curve approach each other, creating a barrier against strand crossing that allows numerical minimization.
\end{itemize}

A well-defined knot energy typically has the effect of preventing two distinct parts of a knot from getting either too close or too far from each other. Many knot energies have been proposed \cite{buck1997energy,von2017dynamics, yu2021repulsive,lipton2022stationary}, but the most-studied is likely the Möbius energy. O'Hara proposed a family of knot energies based on a inverse-power repulsive potential between any two points on a curve balanced by a term depending on the distance between the points along the curve \cite{o1994energy}. When this power is -1 we recover something analogous to the Coulomb energy, which does not provide a barrier to strand crossing. While the Coulomb energy is physically motivated, it is minimized by a charged circle of infinite radius which makes it a bad knot energy. Freedman, He, and Wang pointed out that when the power of the O'Hara energy is -2, the energy is invariant under conformal and Möbius transformations, terming this the Möbius energy \cite{freedman1994mobius}. It is defined within a single closed curve $\gamma$ as follows:

\begin{equation}
E_{s}=\int_{\gamma(u)}\int_{\gamma(v)}\left(\frac{1}{\left(\gamma(u)-\gamma(v)\right)^2}-\frac{1}{D_\gamma(u,v)^2}\right)
|\dot{\gamma}(u)||\dot{\gamma}(v)|du\ dv
\end{equation}
$D_{\gamma}(u,v)$ is the distance along the curve between $\gamma(u)$ and $\gamma(v)$, and the overdots represent tangent vectors. The first reciprocal term that depends on the absolute distance can be thought of as a repulsive interaction that prevents distinct segments of the curve from being too close. The term that depends on the distance along the curve can be thought of mathematically as a regularizer that prevents the integral from diverging, or physically as a pseudo-tension that favors straight segments between points of fixed distance.

When multiple curves are connected in a link, the above expression is termed the Möbius self-energy and the cross-energy can be calculated between any pair of curves:

\begin{equation}
    E_M=\sum_i{E_{s,i}}+\sum_{i\neq j}{E_{c,ij}}=\sum_i{E_{s,i}}+\sum_{i\neq j}\int_{\gamma_i(u)}\int_{\gamma_j(v)}\left(\frac{|\dot{\gamma_i}(u)||\dot{\gamma_j}(v)|du\ dv}{\left(\gamma_i(u)-\gamma_j(v)\right)^2}\right)
\end{equation}
There are several conventions in how the Möbius self- and cross-energies are defined (see footnote \footnote{Some authors, including Freedman-He-Wang \cite{freedman1994mobius} and Rawdon \& Simon, define the energy such that the circle takes a value of 4. Others, such as Kusner and Sullivan \cite{kusner1998mobius}, define it such that circle has zero energy and every other curve has an ``excess energy'' relative to that minimum. When discussing the cross-energy of links, one convention, used by Freedman-He-Wang and Agol-Marques-Neves \cite{agol2016min}, counts each pair of line elements once such that a Hopf link has a minimum cross energy of 2$\pi^2$. Another, used by Kusner and Sullivan, counts each pair of line elements twice, as in the self-energy, such that the Hopf link has a minimum cross energy of 4$\pi^2$. We work with the 4/$4\pi^2$ convention, for comparisons to other works. Arguments can be made for either convention: if each pair of line elements is always double-counted, the Hopf link should have a cross-energy of $4\pi^2$. If we begin with a definition for the cross-energy of two curves that avoids double counting, and then translate the two curves to perfectly overlap (as is done when using the Gauss linking integral to describe the space writhe or average crossing number), the self-energy is recovered without double-counting.}), we work with the conventions that the Möbius energy of a circle is 4 and the minimum cross-energy of the Hopf link is $4\pi^2$. Since its definition and exploration, many properties of the Möbius energy and its minimizers have been proven \cite{blatt2022hara,blatt2012boundedness, diao2012properties}, and several projects have been undertaken to find the minimal configurations of certain knots \cite{kim1993torus,kusner1998mobius, lipton2022complex}. Because of the conformal invariance of the Möbius energy, a given knot does not have a unique minimizer, but a family of embeddings that minimize the energy.

The energy of a knot is typically minimized with a gradient descent algorithm based on a piecewise-linear discretization of a knot into Cartesian coordinates, which raises questions about the best to define and measure the discretized energy and whether it converges to the smooth value. Avoiding these questions, Simon defined the Minimum Distance (MD) energy \cite{simon1994energy}, based on the sum of the inverse squares of the minimum distance between non-adjacent line segments of a closed polygonal curve:
\begin{equation}
    E_{MD}=\sum_{i,\ j,\ |i-j|>1}\frac{\ell_i\ell_j}{(MD_{ij})^2}
\end{equation}
where $\ell_i$ and $\ell_j$ are the length of line segments i and j, and $MD_{ij}$ the minimum distance between the two line segments. The smallest polygon that has a defined MD energy is the quadrilateral, and the energy is minimized by a square with an energy of 4 (factor-of-two conventions exist here as well). An interesting feature of the MD energy is that there is an optimal number of edges that will minimize the MD energy for a given knot. While the unknot is minimized by a square, it is known that the hexagonal trefoil has a greater minimum energy than the heptagonal trefoil and that the likely minimizer is the 14-gon \cite{simon1994energy}. There has not been an extensive exploration of the minimum MD energies of various knots, with a few detailed investigations of regular polygons and hexagonal trefoils. Links have not typically been discussed in the context of MD energy, but since two edges on different polygons cannot be adjacent, no modification of the energy definition is required.

While one of the goals of a knot energy is to allow ``nice'' embeddings, it is not always clear what those nice embeddings are even when the energy is minimized.   A minimum energy can be proven without prescribing how a curve can achieve it. If handed a length of wire and asked to bend it into a nice knot, one cannot invert the space containing the knot or embed it in four dimensions before starting to bend. Conversely, an algorithm can reach a minimum but not provide insight as to why that configuration is minimal. The goal of this manuscript is to bridge that gap and provide geometric guidance on how minimal link energies can be achieved. We focus on links of multiple planar unknots such as the Hopf link and Borromean rings, building up in complexity towards extended chains and planar chainmails, for reasons discussed in the next paragraph. We describe embeddings of these knots with respect to free geometric parameters and examine which values of these parameters will, at least locally, minimize the energy of the link. We will first discuss considerations and implications of the Möbius energy for smooth curves, then of the MD energy for polygonal curves, before discussing efficient chainmail constructions.

I will speak briefly on the inspiration for this work. An oft-stated motivation for knot energies is to understand knots in DNA molecules, but as a researcher in the field for the past decade I have not found this to be the case. DNA molecules minimize free energy rather than just energy, and their ``long-range'' interactions are governed by a quickly-decaying Yukawa potential which does not have an associated knot energy. However, in my investigation of kinetoplast DNA (molecular chainmail) I found that ropelength-minimizing Hopf linked networks showed the same Gaussian curvature as polymer networks with the same lattice topology \cite{klotz2024chirality}. This motivated me to investigate Möbius energy-minimizing chainmail networks, after which I realized that much simpler questions like ``how far apart should the two rings in a link be'' did not have an answer, leading to this manuscript.

\subsection{Methods and Algorithms}

To compute the Möbius energy, curves are discretized into polygons of typically a few hundred points. Each point acts as a ``charge'' and the contribution to the Möbius energy of each pair of charges is computed during a nested \textit{for} loop. The distance along the curve between two points is simply calculated from the sum of the distances between each pair of adjacent vertices between the two points. If this is greater than half the total length of the curve, it is subtracted from the total length.

For the Möbius energy of a curve, our algorithm has the same accuracy with respect to the number of vertices as that reported by Kim and Kusner \cite{kim1993torus}, with an underprediction that decreases with increasing vertex count. With 360 vertices, the energy of a circle is found to be 3.9607, with 720 it is 3.9804. The cross-energy calculation is more accurate because computing the distance along the curve is not required. A Hopf link with 180 vertices in each circle yields an energy that is 0.999897 of the known minimum. 

To compute the MD energy, we require the minimum distance between two line segments. This is found using an algorithm from Dan Sunday's now-defunct C algorithms website that was transcribed into MATLAB by a user named only Nick. The original description of the algorithm can be found on the Internet Archive \cite{sunday}.

Our results are typically based on describing a given link as a function of a few free parameters, instantiating those links with a discrete number of vertices, then minimizing the energies with respect to the free parameters. For minimization we use two native MATLAB functions: fminsearch, which minimizes a multi-parameter function based on the Nelder-Mead Simplex Method \cite{lagarias}, and fminbnd, which minimizes a single-parameter function between two bounds using a golden-section search \cite{brent2013algorithms}. We have included a sample MATLAB script as ancillary data. Typically when we achieve a minimal Möbius energy, we re-calculate it with double the vertex count to evaluate convergence. The minimizing parameters, e.g. sizes and positions of components, are typically insensitive to vertex count.

As a preliminary investigation we wrote a LAMMPS script that reduces the Möbius energy of a knot using a long-range inverse-square repulsive energy but relying on inextensibility rather than the contour term. While not as direct a path  to minimization as a dedicated Möbius minimizer, this may serve as a more standard tool for users interested in carrying out future investigations. As a proof of concept, a 100 vertex trefoil was reduced to a Möbius energy of 74.88 (after spline-interpolation to 1000 vertices), compared to the reside minimum of 74.41 \cite{kim1993torus}. 

\section{Möbius Energy Minimization}

\subsection{Hopf Links}

\begin{figure}
    \centering
    \includegraphics[width=0.8\linewidth]{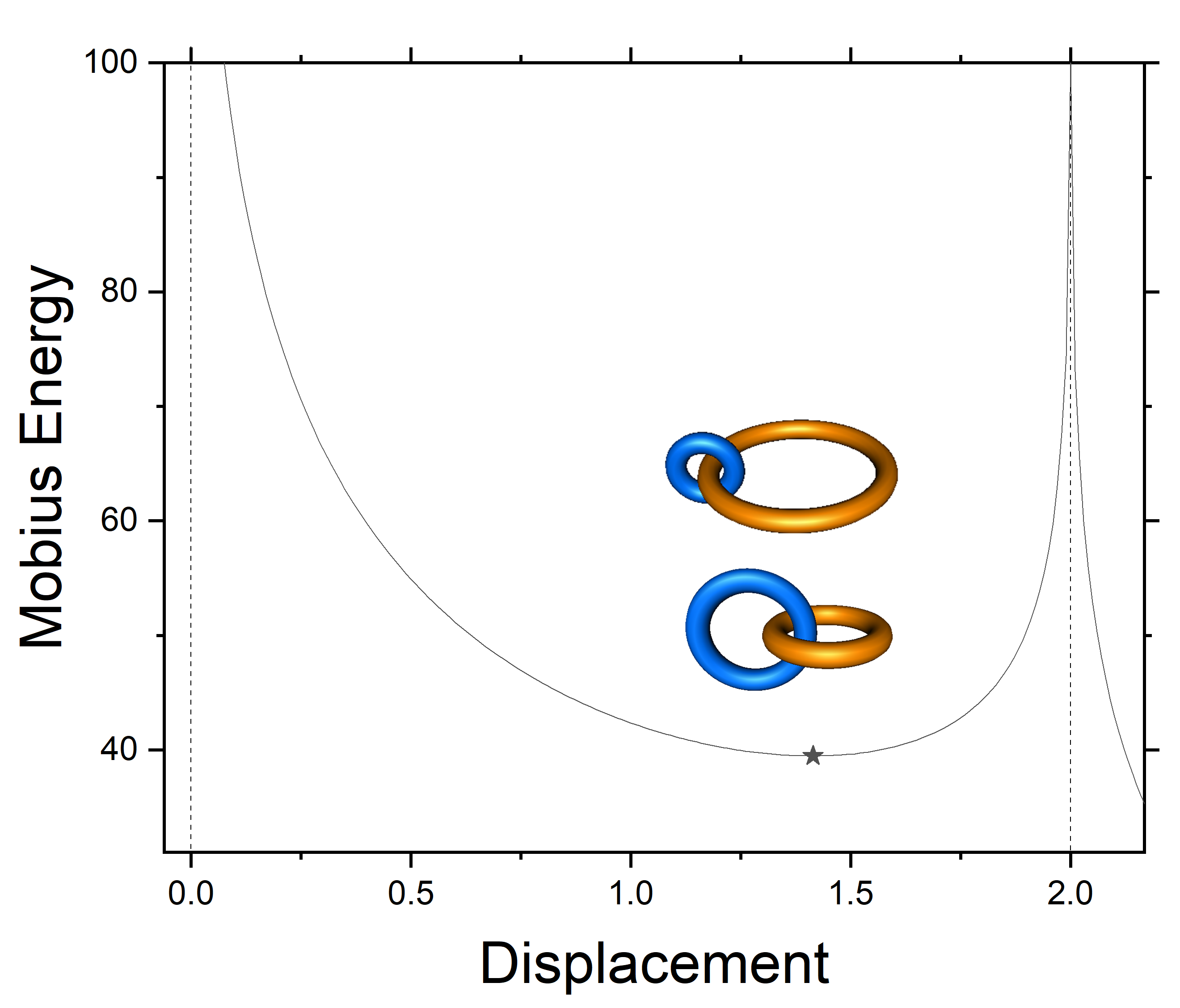}
    \caption{Möbius energy of two unit circles with perpendicular normals as a function of the distance between their centers along a mutually perpendicular axis. The energy is minimized at a displacement of $\sqrt{2}$ with a value of $4\pi^2$. The energy diverges when the two circles overlap. Two energy minimizers are shown, the bottom with congruent circles corresponding to the separation in the plot, the top with an asymmetric radius ratio.}
    \label{fig:hopf}
\end{figure}

It was initially conjectured that the minimum cross-energy of a Hopf link is $4\pi^2$ (or $2\pi^2$), and this was proven in 2014 \cite{agol2016min}. The initial conjecture was based on the ability to use a Möbius transformation to transform the entire z axis into a circle, and the proof based on two circles each lying in two-dimensional planes of four-dimensional space. Neither scenario arrives at the specific coordinates of a Möbius-minimizing Hopf link. If one computes the energy of two unit circles with perpendicular normals that pass through each others' centers, the energy will be 7.3\% above the proven minimum. To find the coordinates in three-dimensional space that minimize the energy of the Hopf link, we explicitly solve the integral for the Möbius cross-energy of two circles with perpendicular normals, one of which is displaced from the origin by a distance $\Delta$ along the cross product of the two circles' normal vector.   
We can parameterize a Hopf link by a unit circle $\vec{r_1}$ and another circle $\vec{r_2}$ with arbitrary radius $\alpha$:
\begin{equation}
    \vec{r_1}(\theta)=\left\langle\cos(\theta),\sin(\theta),0\right\rangle,\ \ \ \ \vec{r_2}(\phi)=\left\langle\alpha\cos(\phi)+\Delta,0,\alpha\sin(\phi)\right\rangle
\end{equation}
In the simplest case, $\alpha=1$, and $0<\Delta<2$. The cross-energy of the link in this case can be solved: 
\begin{equation}
    E_c=2\int_0^{2\pi}\int_0^{2\pi}\frac{d\theta d\phi}{\left(\cos\theta+\Delta-\cos\phi\right)^2+\sin^2\phi+\sin^2\theta}=\frac{-16\pi}{\Delta^2-4}E_K\left(-8\frac{\Delta^2-2}{\left(\Delta^2-4\right)^2}\right)\approx4\pi^2+2\pi^2(\Delta-\sqrt{2})^2
\end{equation}
Here, $E_K$ is the complete elliptic integral of the first kind. The solution to the double integral was determined by the sophisticated method of massaging partially-helpful \textit{Mathematica} outputs combined with pattern recognition. The energy as a function of separation is plot in Fig. \ref{fig:hopf}. The expression in Eq. 6 achieves the known minimum energy when the two unit circles are separated by a distance of $\sqrt{2}$. This is an explicit example of a link for which the ropelength-minimizing embedding, which has the circles separated by a distance of 1, does not minimize the Möbius energy. The two linked circles also need not have the same radius to minimize energy, the Möbius cross-energy between a unit circle and one with radius $\alpha$ is minimized to the same $4\pi^2$ when $\Delta=\sqrt{1+\alpha^2}$, a relation that was determined by pattern recognition and verified by \textit{Mathematica}.

\begin{equation}
    E_c=2\int_0^{2\pi}\int_0^{2\pi}\frac{\alpha\ d\theta d\phi}{\left(\alpha\cos\theta+\sqrt{1+\alpha^2}-\cos\phi\right)^2+\sin^2\phi+\alpha^2\sin^2\theta}=4\pi^2
\end{equation}

The energy is minimal with respect to unsurprising transformations, e.g. changing the perpendicular inclination of the two circles, or translating them along any directions besides that of the cross product of their two normal vectors.

\subsection{Hopf Tambourines and Ropelength}

Three circles of equal radius can be linked in a chain, with the two outer circles each centered a distance of 1.58 from the inner circle, and the energy will be locally minimal at 102.4. This is 12.6\% above the minimum energy of $12+8\pi^2$ because the cross-energy of the two outer circle contributes non-negligibly to the total energy. If, however, the two outer circles are much smaller, they would each have the same Möbius self-energy and cross-energy with the larger central circle, but effectively zero cross-energy with each other. This would give the entire link close to the minimum possible energy.

\begin{figure}
    \centering
    \includegraphics[width=0.6\linewidth]{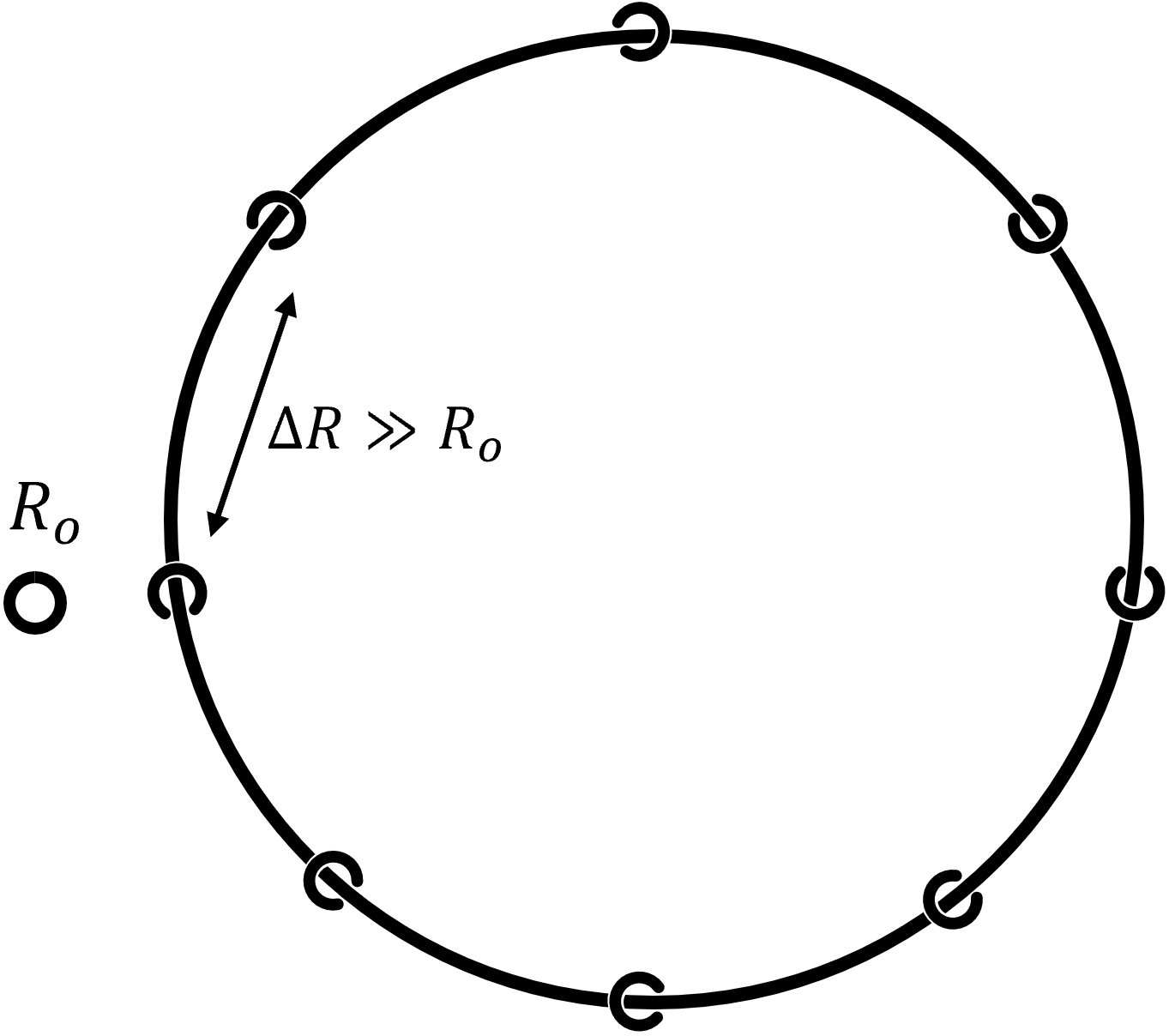}
    \caption{A ``tambourine'' configuration of many small circle Hopf linked to a larger central circle. As the radius of the small circles approaches zero, the energy between each small circle and the central one stays constant at $4\pi^2$ but the energy due to the interactions between small circles approaches zero. The tambourine is conjectured to be the minimum Möbius energy for any link with an even crossing number.}
    \label{fig:tambo}
\end{figure}

If more small circles are linked to the large central circle, they will all be separated by far enough relative to their radius to have negligible cross energy. We will refer to this configuration as a ``tambourine,'' an image evoked by Fig. \ref{fig:tambo}. For two circles separated by a distance $R$ much greater than their radius $r$, they can be treated as ``point charges'' with cross-energy $8\pi^2r^2/R^2$ (it is somewhat of a coincidence that this formula produces the correct minimum when $R/r=\sqrt{2}$). If a total of $N$ small circles are linked to a central large one, the total Möbius energy is:

\begin{equation}
    E_M>4(N+1)+4\pi^2N
\end{equation}

The crossing number is $C=2N$ which allows the energy to be written, when $C$ is large, as:

\begin{equation}
    E_M\approx(2\pi^2+2)C
\end{equation}
The ropelength of a knot is the minimum ratio of the contour length to the thickness of a knot, where the thickness is defined based on the radius of the smallest circle that can pass through three points on a knot. Rawdon and Simon \cite{rawdon2002mobius} showed that the ropelength of a knot ($L$) is constrained to its Möbius energy:
\begin{equation}
    L>\left(\frac{E_M}{4.57}\right)^{3/4}
\end{equation}
The denominator has an exact expression, and reduces to 3.63 when $C$ is large. If the Hopf tambourine has the minimum Möbius energy among every knot or link with $C$ crossings, this constrains the universal lower bound on ropelength:
\begin{equation}
    L>\left(\frac{2\pi^2+2}{4.57}C\right)^{3/4}\approx3.22C^{3/4}
\end{equation}
This would be a significant improvement over the current lower bound prefactor of 1.10. When $C$ is large, this would become 3.82. For odd crossing numbers, one of the linked circles can be replaced by a trefoil, increasing the energy by more than 70, which will not be significant when $C$ is large. This being a true lower bound for ropelength relies on a few assumptions:
\begin{itemize}
    \item \textbf{Rawdon and Simon's bound applies to links.} Links are not explicitly discussed in that work, but because Hopf chains have a ropelength linear in crossing number, the 4/3 bound between ropelength and energy is not violated by them.
    \item \textbf{The $4\pi^2$ convention for cross-energy applies to Rawdon and Simon's argument.} If the $2\pi^2$ Hopf convention is used, the coefficients are reduced to 2 and 2.43, still an improvement over 1.10. However, a lower bound cannot rely on convention choice.
    \item \textbf{No other knot or link with $C$ crossings has a lower energy than a Hopf tambourine.} $2\pi^2+2$ is approximately 21.7, and the most efficient knots and links that have been studied have energies per crossing of around 24. Although these energies have not been extensively studied beyond 8 crossings, it is conceivable that a more efficient knot could exist. If the $2\pi^2$ convention is correct, this ropelength bound is more likely to be true. Because non-alternating knots typically have lower ropelength than alternating knots with the same crossing number, one might expect them to have lower Möbius energies. However, if there were a family of knots whose minimal Möbius energy scaled sub-linearly, then the ropelengths of those knots would grow with a power that is below 3/4 which is known to be impossible.
\end{itemize}

Currently, a stronger lower bound on ropelength from Diao \cite{diaorope} exists for knots below 1850 crossings. If the tambourine geometry constrains ropelength as described, it would only be improved upon by Diao's bound below 48 or 220 crossings (depending on whether the prefactor is 3.22 or 2), rather than 1850.

\subsection{$6_3^3$, Borromean Rings, and Others}

\begin{figure}
    \centering
    \includegraphics[width=1\linewidth]{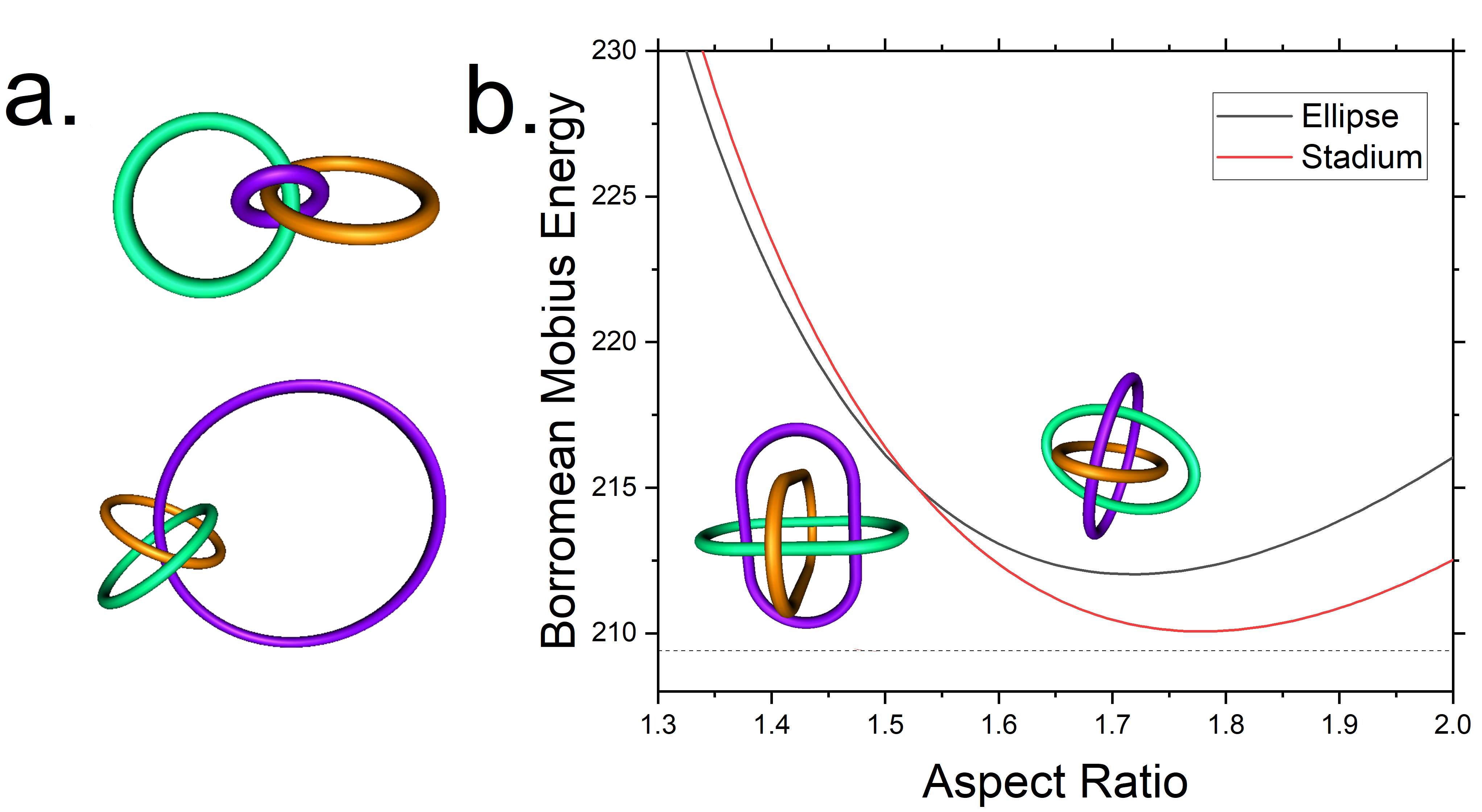}
    \caption{a. Minimal configurations of the $6_3^3$ link. The top figure, with the larger unit circles separated by a distance of $\sqrt{3}$ and inclined at 60 degrees relative to each other, and the smaller circle with half the radius and inclination, achieves the known minimum. The bottom configuration only achieves the minimum as the larger circle's radius tends to infinity.  b. Energy of Borromean rings consisting of three perpendicular ellipses or stadium curves as a function of their aspect ratios, with representative images of the local minimizers. The dashed line is the numerical minimum found by Kusner and Sullivan \cite{kusner1998mobius}.}
    \label{fig:borro}
\end{figure}




There are other nontrivial links that can be constructed from planar unknotted curves and optimized. The non-alternating $6_3^3$, also T(3,3), has a known energy minimum found by Kim, Strengle and Kusner at $12+8\sqrt{3}\pi^2\approx148.75$ \cite{kim1993torus}. Kusner and Sullivan reached 148.8 with their discretized minimizer. By initializing the link as two unit circles with a variable separation and inclination, and a third circle at their midpoint with half the inclination and a variable radius, we find that the minimum energy is achieved when: the two unit circles are separated by a distance of $\sqrt{3}$, they are inclined at 60 degrees about their axis of separation relative to each other, the small circle has a radius of 1/2, and an inclination of 30 degrees relative to both. The circles can also be initialized such that the solver is unstable against runaway growth of one of the components. The energy of this configuration approaches the minimum as the radius tends to infinity, at which point the inclination angles reach those seen in the stable configuration. This is related to the non-existence of minimizers for composite knots \cite{kusner1998mobius}. The stable and unstable configurations of $6_3^3$ are seen in Fig. \ref{fig:borro}a.

Borromean rings cannot be formed from three circles \cite{lindstrom1991borromean}, but can be formed from ellipses are other elongated shapes.  Three congruent ellipses with mutually orthogonal normal vectors and common centers will form Borromean rings, and such a configuration is minimized when the ellipses have an aspect ratio of 1.71. Such a configuration is minimal with respect to small translations and rotations of the ellipses. However, a stadium curve has a lower Möbius energy than an ellipse with the same height and width beyond an aspect ratio of 1.74, and three Borromean stadia with an aspect ratio of 1.78 will have an energy that is 0.5\% lower than the best ellipse configuration. Shapes such as stretched squircles and rounded rectangles were found to have lower energy than ellipses, but greater than stadia. The minimal parameters are slightly different if the $2\pi^2$ conventions for the cross-energy is used, lowering the optimal aspect ratios. The minimal energy of the stadia Borromean rings is slightly higher (210.1 vs 209.4) than the best value determined by gradient optimization. A picture provided by Gunn and Sullivan \cite{gunn2008borromean} suggests that the ideal shape is a stadium with slightly rounded sides. The best reported Möbius cross-energy of Borromean rings is within a factor of 1.00004 of $20\pi^2$ or $10\pi^2$ in the lesser convention, which may warrant investigation.

Attempts were made to minimize the $4_1^2$ (Solomon/T(4,2)), $5_1^2$ (Whitehead), $6_1^3$, and $8_3^4$ (a non-alternating Hopf loop), but none of the local minima were close enough to the best numeric results to provide insight. Some of these require non-planar components to be embedded, which are typically parameterized by higher Fourier modes in the z-direction.

Torus knots parameterized by the standard harmonic formula have an energy primarily determined by the major and minor radius of the torus. We can minimize torus knots with respect to these parameters, and find values slightly larger than the minima determined by Kusner and Strengle's residue formula (for example, 80.08 for the trefoil instead of 74.41). The configurations are typically minimal with respect to anisotropy of the curves (e.g. scaling in the z direction), but the energy can be reduced with higher harmonic modes.

\section{Minimum Distance Energy}

\subsection{Hopf Links and Borromean Rings}

We can find exact parameters that locally minimize the MD energies of simple links. The minimum distance energy of a Hopf link consisting of two squares can be computed from the sixteen pairwise interactions of the four line segments they each consist of, plus the energy of each square. Likewise, the MD energy of Borromean rings is thrice the pairwise energy of two interpenetrating rectangles, plus thrice the energy of each rectangle. The Hopf energy can be minimized with respect to the distance between the two squares and the Borromean energy by the aspect ratio of the rectangles. Both are minimized when the two quadrilaterals have perpendicular normals and share common centers in at least two dimensions. In both cases, while the parameters for energy minimization can be solved exactly, the MD energy of the link is \textit{not} minimized by quadrilaterals, but by pentagons for the Hopf link and decagons for Borromean rings.

\begin{figure}
    \centering
    \includegraphics[width=1\linewidth]{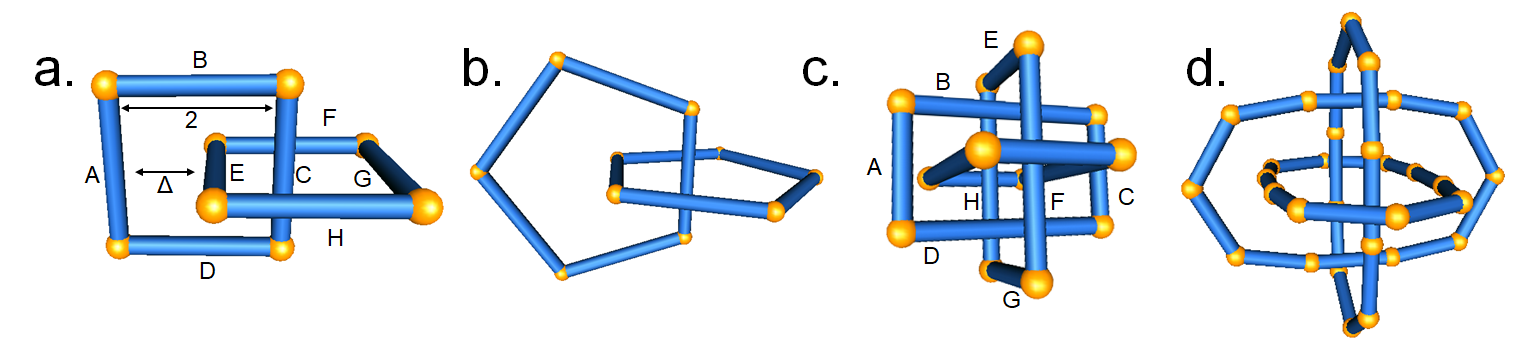}
    \caption{a. A Hopf link consisting of two squares separated by an MD energy minimizing distance $\Delta$ of 1.203. The labeled line segments are used in the text to derive that distance. b. A Hopf link of two pentagons with slightly lower MD energy. c. Borromean rings consisting of three rectangles with an MD-minimizing aspect ratio of 1.756, labelled similarly to a. d. Borromean rings constructed of decagons with optimized vertex positions, our best minimizer.}
    \label{fig:mdfig}
\end{figure}

Figure \ref{fig:mdfig}a shows a square with sides of length 2 labeled $A, B, C, D$ linked to one labeled $E,F,G,H$. Square $ABCD$ lies in the XY plane and is centered on the origin, and $EFGH$ lies in the XZ plane and its center is displaced along the X direction by a distance $\Delta$. Within each square, energy arises from non-adjacent segment pairs $AC$, $BD$, $EG$, and $FH$, all of which have minimum distance 2. We must find the minimum distances between each segment on ABCD and those on EFGH. From symmetry, the minimum distances of $CF, CH, BE,$ and $DE$ are all 1. Similarly, $AF$, $AH$, $BG$, and $DG$ are the same at $\sqrt{1+\Delta^2}$ and $BF, BH, DF,$ and $DH$ each have a minimum distance of $\sqrt{2}$. $CG$ and $AE$ have a minimum distance of $\Delta$. Finally, $AG$ has a minimum distance of $2+\Delta$ and $CE$ has a minimum distance of $2-\Delta$. The total cross-energy of the two squares is:

\begin{equation}   E_H=E_{AG}+E_{CE}+2E_{AE}+4E_{BE}+4E_{BF}+4E_{AF}=4\cdot2\cdot\left(\frac{1}{(2+\Delta)^2}+\frac{1}{(2-\Delta)^2}+\frac{2}{\Delta^2}+\frac{4}{1^2}+\frac{4}{2}+\frac{4}{1+\Delta^2}  \right)
\end{equation}
The factor of 4 in front arises from the product of the side length in the energy definition, the factor of two is from the pair-counting convention. This has an exact minimization of approximately 1.203, such that $x=\Delta^2$ is the real positive solution of the quintic $2x^5-10x^4+73x^3-48x^2-40x-32$. The minimal energy is approximately 93.5, of which 8 is the self-energy of the squares. This is minimal with respect to translation, rotation, and elongation of both squares and rescaling of their relative sizes. With respect to a size asymmetry, the energy is qualitatively hyperbolic and asymptotically linear. This configuration of two squares is not, however, the minimal-energy Hopf link, as a configuration of two linked regular pentagons has a slightly lower energy 90.93 (Fig. \ref{fig:mdfig}b). 


Mutually perpendicular squares cannot form Borromean rings, but rectangles can. We can compute the pairwise cross-energy of two rectangles, add it to the self-energy of one rectangle, and triple the result for the total cross-energy For two of the interpenetrating rectangles making up Borromean rings, we can align one with its short side along the X axis and its long side along the Y axis, and the other with its long side along the X axis and its short side along the Z axis, both centered on the origin. If the rectangles have side lengths A and B and aspect ratio $\alpha=B/A$ the MD energy of each rectangle is:
\begin{equation}
    E_{R}=2(\frac{A^2}{B^2}+\frac{B^2}{A^2})=2(\frac{1}{\alpha^2}+\alpha^2).
\end{equation}
Defining the four sides of each rectangle according to Fig. \ref{fig:mdfig}c, with short length 1 and long length $\alpha$, we have five distinct distances each describing four pairs of sides:
\begin{align*}
    AE=AG=CE=CG=\sqrt{(\frac{\alpha}{2})^2+(\frac{\alpha-1}{2})^2}\\
    AH=AF=CH=CF=\sqrt{(\frac{\alpha}{2})^2+(\frac{1}{2})^2}\\
    BH=BF=DH=DF=\frac{1}{2}\\
    BG=DE=\frac{\alpha-1}{2}\\
    BE=DG=\frac{\alpha}{2}+\frac{1}{2}
\end{align*}
After some simplification the total energy can written as:
\begin{equation}
    E_B=6\cdot\left(17\alpha^2+\frac{16\alpha}{\alpha^2+1}+\frac{16}{2\alpha^2-2\alpha+1}+\frac{8\alpha}{(\alpha-1)^2}+\frac{8\alpha}{(\alpha+1)^2}+\frac{1}{\alpha^2}\right)
\end{equation}
The value of $\alpha$ that minimizes this energy is the real positive solution of an 18th-order polynomial, and is approximately 1.756, comparable to the aspect ratios of Möbius minimizers. This is quite close to the aspect ratio of the stadium curve that minimizes the Möbius energy of Borromean rings. The minimal energy is 542.6.

Because an answer arising from the root of a high-order polynomial is not particularly more insightful than a number determined from numerical minimization, we proceed numerically. Like the Hopf link, the quadrilateral does not minimize the MD energy for all Borromean links. Using regular n-gons that have been stretched in one dimension and rotated to have favorable orientation (we can treat the angle as a free parameter but it is typically $\pi/n$ relative to their harmonic parameterization), we find that decagons have lowest energy minimum at 406, with odd-gons having slightly higher energy than even-gons. The ideal aspect ratio of the decagon is 1.74, but as the number of sides increases the ratio approaches roughly 1.715, the same as the best ratio for the Möbius energy of Borromean ellipses. This is consistent with the known result that an MD-minimizing polygon approaches a Möbius-minimizing curve as its vertices increase in number \cite{rawdon2010error}. Since the decagons are planar and have four-fold symmetry, the energy can be minimized as a function of only four variables. Specifically, the corners of the decagon may be thought of as numbers on a clock missing 3 and 9, and 12 o'clock may be placed at (0,1) without loss of generality. After minimization the ideal coordinates of 1 and 2 o'clock are (0.44, 0.75) and (0.52, 0.27), and the rest of the vertices may be placed by symmetry. The energy is 387.9 and the shape is almost hexagonal (Fig. \ref{fig:mdfig}d).

The author is generally skeptical of the use of generative artificial intelligence, but on a whim asked ChatGPT ``Describe the two polygons that minimize the minimum distance energy of the Hopf link.'' To my surprise, ChatGPT described regular pentagons as the minimizing shape, but incorrectly mentioned that they passed through each others' centers. It cited a paper on ropelength by Cantarella et al. \cite{cantarella2002minimum}, which notably does not discuss MD energy or make this claim. While ChatGPT's arrival at pentagons may be an example of the stopped-clock principle, it may be possible that I have missed relevant literature. When asked about Borromean rings, ChatGPT incorrectly stated that regular hexagons minimize the energy, citing erroneously the same paper.



\subsection{Linear Hopf Chains}

A linear chain of Hopf linked unknots has a known ropelength, but its knot energies have not previously been considered. Here we examine the energy of chains of odd numbers of rectangles. The simplest such case has three components. If the components are congruent, their energy is minimized by a circle or square centered on the origin in the XY plane, and two others in the XZ plane translated an equal distance in the X direction. For the Möbius energy, the ideal spacing is 1.58 and the energy is 12.8 \% above the absolute minimum. For the MD energy, the distance between the squares is similar, approximately 1.57. As more congruent components are added, if the periodicity between successive components is held constant then the spacing slowly increases with the number of components. Congruent shapes however do not minimize energy.

Because of Möbius energy of a Hopf link of two circles is independent of the relative size the circles, the energy can be minimized by a pathological configuration of a central unit circle, two much smaller circles linked on either side such that their cross energy is effectively zero. A five component link can be constructed by linking even smaller circles to the already-small circles linked to the central unit circle, etc. The situation is more interesting when considering the MD energy.

\begin{figure}
    \centering
    \includegraphics[width=1\linewidth]{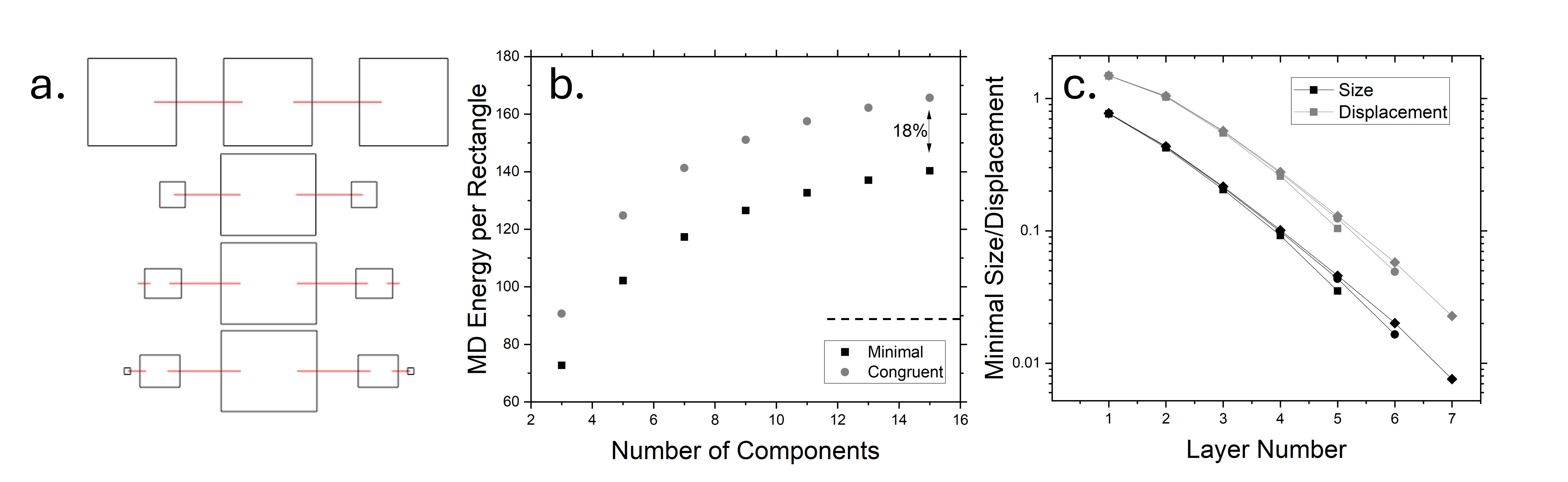}
    \caption{a. Top image shows the Minimum Distance energy minimizing configuration of a 5-square Hopf chain if the size and spacing of the squares is fixed. Below images show energy minimizing 5, 7, and 9 rectangle Hopf chains when the size, aspect ratio, and spacing of each rectangle is allowed to vary. b. The MD energy per rectangle as a function of a number of rectangles, with and without a congruency constraint. The dashed line shows the minimum possible energy. c. The size and displacement from their neighbors of each pair of rectangles in increasing layers from the center. The three curves of each set show the minimum for X, Y, and 15 total rectangles. }
    \label{fig:squarechain}
\end{figure}

The top image in Figure \ref{fig:squarechain}a shows the minimal MD energy configuration of a congruent five-square Hopf chain, with the red squares lying transverse to the page. Assuming that the MD energy is minimized by by an extended configuration of perpendicular rectangles, a Hopf chain of $2N+1$ components can be described by $3N+1$ free parameters: each symmetric layer of rectangles has a distinct area, a distinct displacement from their inner neighbors, and a distinct aspect ratio. The central rectangle can be treated as having unit area and centered on the origin, but has a free aspect ratio. The remaining images in Fig. \ref{fig:squarechain}a show the minimal configurations with these degrees of freedom for 5, 7, and 9 components. The 5-component link is considerably smaller than the congruent one, and each successive layer of rectangles is considerably smaller than the previous. The energy per rectangle approaches a plateau of about 140 as each new pair of rectangles is far enough apart that they add minimal excess cross energy. The congruent links have about 18\% more energy than the minimized links, which are about 57\% greater than the minimum (of 4 per square plus 85.5 per linkage). 

While not as drastic as Möbius minimizers, the size and displacement of each additional layer of rectangles shrink exponentially (Fig. \ref{fig:squarechain}c), each about a factor of 2-3 smaller than the previous. While there are small changes of the interior rectangles when new layers are added, the sizes and displacements are largely self-similar. The aspect ratio of every rectangle in these chains was approximately 1.2 (elongated on the x-axis), except those on the outermost layer which were slightly compressed at 0.97. Interestingly, this exponential decrease in size and displacement implies that as the number of components in the chain tends to infinity, the width subtended by the entire chain will remain finite. The width of the 15-component link is approximately 7.2, while the width of the 9 component link is approximately 6.5. This makes these energy minimizing links qualitatively similar to wild knots. Although we have only demonstrated this with rectangles, we suspect similar behavior from other polygons. We do not know if a chain of pentagons would be the global minimizer.

\section{Energy-Minimizing Chainmail Networks}

Chainmails can be thought of as two dimensional networks of linked unknots, which may be alternating or non-alternating. Alternating chainmail networks were recently proven to be L-space links \cite{agol2023chainmail}, and the relationship between the alternatingness of their network topology and the Gaussian curvature of their embeddings has been considered in the context of kinetoplast DNA \cite{klotz2024chirality}. Globally minimizing the energy of a chainmail network may be difficult due to the large number of parameters and the computational challenge presented by a large network. Here we consider the Möbius and MD energies of two types of chainmail networks, each minimized with respect to two free parameters. European 4-in-1 style chainmail consists of a square lattice of rings, each linked with its four nearest neighbors in an over-under-under-over-under-over-over-under pattern (other sequences of over and under do not produce a flat network). Square Japanese-style chainmail (also known as kusari katabira) consists of a square lattice of planar rings, each linked to its nearest neighbors by a ring oriented in a different plane. Every ring in 4-in-1 chainmail is the same, whereas rings in Japanese chainmail will have a valence of either 4 or 2.

4-in-1 chainmail is constructed by unit circles at the vertices of a square lattice in the XY plane with periodicity $D4<2$, and rotating the circles in each odd row about the y-axis by angle $\theta_4$, and the circles in each even row by $-\theta_4$. After initializing a network with $N^2$ total links, the energies can be minimized with respect to \textit{D4} and $\theta$. Japanese chainmail is constructed by tiling circles in a square lattice in the XY plane by $DJ>2$, and taking circles in the YZ or XZ planes and placing them at the midpoints between the planar circles. The energies can be minimized with respect to the spacing between the planar circles and the size of the linking circles ($LJ$). For minimizing the MD energy, the circles in both cases can be replaced with squares. The minimal configurations reached are only minimal with respect to the two free parameters, adding a third parameter corresponding e.g. to the aspect ratio of the components may allow further reduction, and insights from examining Hopf chains indicate that a local minimum may not be significantly more energetic than a global minimum.

\begin{figure}
    \centering
    \includegraphics[width=1\linewidth]{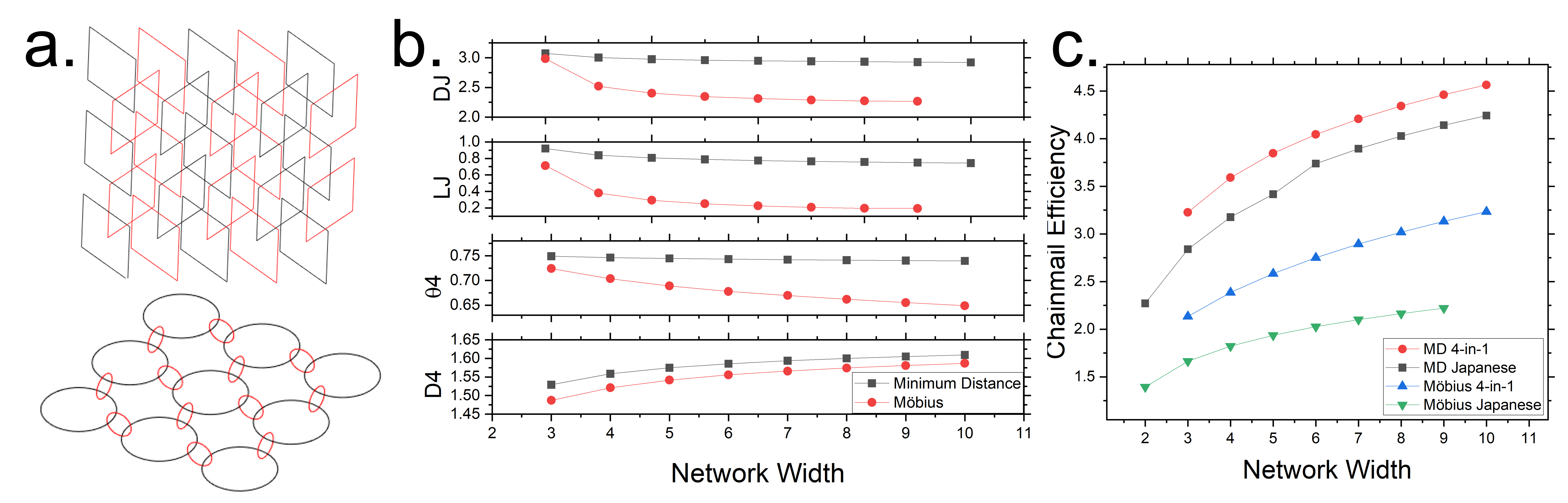}
    \caption{a. Top shows a 5x5 network of 4-in-1 chainmail minimizing the minimum distance energy. Bottom shows a 3x3 network of Japanese chainmail minimizing the Möbius energy. b. Minimal parameters for chainmail energies as a function of the network size. The top two plots show the lattice constants and linker size of Japanese chainmail. The bottom two plots show the tilt angle and neighbor distance for 4-in-1 chainmail. Black corresponds to minimum distance energy and red corresponds to Möbius energy. c. Efficiency, defined as the ratio of the minimized energy to the absolute minimum, for the four types of networks as a function of network size.}
    \label{fig:chainmail}
\end{figure}

Figure \ref{fig:chainmail}a shows a Möbius-minimizing Japanese chainmail network and an MD-minimizing 4-in-1 chainmail network.  The minimizing parameters for each network are shown as a function of the width of the networks in Fig. \ref{fig:chainmail}b, which in all cases varies more for Möbius minimizers than MD minimizers. The distance between planar components in Japanese chainmail (DJ) is approximately 3 for small networks for both energies and very weakly decreases with network size for MD minimizers, while falling towards the minimum of 2 for Möbius minimizers. The corresponding linker size (LJ) decreases similarly with network size. For the 4-in-1 chainmail netowkrs, the angle that each link is rotated with respect to the plane ($\theta_4$) is just below 45 degrees for MD minimizers and small Möbius minimizers, again decreasing very weakly with network size and falling to about 37 degrees for the biggest Möbius minimizers. The lattice constant in both cases rises from about 1.5 to 1.6 over the studied range of sizes. It is an open question as to whether these values converge at an equilibrium for asymptotically large networks, or whether they continue to increase or decrease to their maximum or minimum possible values.

The energy of each network is the sum of the self-energy of each component and the cross-energy of each pair of components. Because we only use minimizing components, we can examine how the total cross energy depends on network size. To do so, we examine the ratio of the total cross-energy to the hypothetical minimum, which is the sum of the minimum energies of all the Hopf linkages in the network. For the Möbius minimizers, this is $4\pi^2$ for each link. For the MD minimizers this is 85.5, the minimum cross-energy of two squares (not two pentagons). The ratio serves as a measure of the excess cross-energy in the network, and is plot in Fig. \ref{fig:chainmail}c. This increases with network size in all cases, sublinearly but not necessarily approaching a plateau. Because the Möbius energy depends less on the relative size of linked components, its excess energy is less than that of MD minimizers. Japanese chainmail networks are generally more efficient than 4-in-1 networks, and this is especially true when considering that they contain nearly three times as many components for a given network width.

One of the motivating questions for this investigation was whether the periodicity and other geometric properties of chainmail networks converges as they become asymptotically large. Does, for example, the linker size of a Möbius-minimizing Japanese chainmail converge on a value near 0.2 as the data suggest, or does it fall towards zero until the planar rings are almost touching? Do the links in a 4-in-1 chainmail get farther and farther apart until they each barely pierce their neighbors? This perhaps may be answered by defining a version of these energies for periodic links. While the sum of inverse squares converges in one dimension, in two dimensions it diverges linearly, which may imply that excess cross-energy must be mitigated by continuously moving rings from each other as the network increases in size.
 
\section{Conclusions}
We have discussed geometric considerations for minimizing the Möbius and Minimum Distance energy of topological links of unknots. We have identified the ideal positions for placing circles, squares, or other polygons to minimize the energies of Hopf links, Borromean rings and other links, as well as extended linear chains and planar chainmail networks. The results we deem most interesting are:
\begin{itemize}
    \item The invariance of the minimum Möbius energy of the Hopf link to the relative sizes of its components implies a way to construct a minimum-energy link for any crossing number, and may imply an improved lower bound for ropelength. 
    \item MD minimizing linear chains subtend a finite length even as the crossing number diverges.
\end{itemize}
Overall we believe these imply that these energies are not ideal for creating ``nice'' configurations of links, even if they are known to produce nice knots. A similar point was made by Yu et al. in their discussion of self-repelling curves \cite{yu2021repulsive}. Other interesting results suggesting further study include the ideal number of vertices in MD-minimizing links, and how these methods might be applied to links of non-planar components or to torus knots. Applications of these findings may include efficient ways to initialize non-interfering chains and chainmails for 3D animation, and possibly to prepare suitable initial configurations for ropelength minimization studies. It is known that the Möbius energy is constrained from below by the absolute Gauss linking number, but links such as the Whitehead link and Borromean rings have nonzero Möbius energy despite zero linking number. Analogous lower bounds could perhaps be derived from Milnor invariants and imply a minimum energy that depends on the degree of ``Brunnian-ness'' of the link. Investigating the $20\pi^2$ Möbius energy of Borromean links may be a suitable first step. 

\section{Dedication and Acknowledgements}
This work is dedicated to my cat Mini (Fig. 7), who passed away while I was working on it. Think of her whenever you minimize a function. I wish to thank Eric Rawdon and Rob Kusner for answering questions about their work on the topic, and apologize to readers for neglected umlauts. The author is supported by NSF grant number 2336744. 

\begin{figure}
    \centering
    \includegraphics[width=0.5\linewidth]{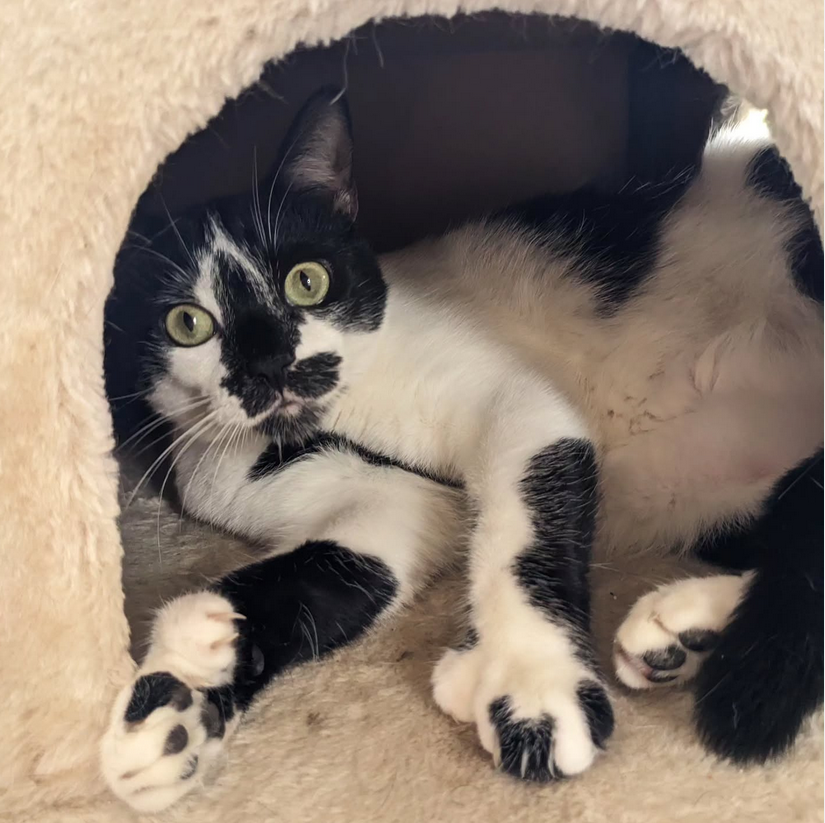}
    \caption{Mini}
\end{figure}

\bibliographystyle{unsrt}
\bibliography{mobrefs}

\end{document}